\documentclass{article}
\usepackage[a4paper]{geometry} % to use the full page
\usepackage{amsfonts,amsmath,enumitem,graphicx,amssymb,mathrsfs,amsthm,color,tikz}
\usepackage[initials]{amsrefs}

\newtheorem{theorem}{Theorem}[section]
\newtheorem{lemma}[theorem]{Lemma}

\newtheorem{corollary}[theorem]{Corollary}

\theoremstyle{definition}
\newtheorem{definition}[theorem]{Definition}

\numberwithin{equation}{section}
\numberwithin{figure}{section}

\renewcommand{\leq}{\leqslant}

\makeatletter
\renewcommand\section{\@startsection {section}{1}{\z@}%
                                   {-3.5ex \@plus -1ex \@minus -.2ex}%
                                   {1.3ex \@plus.2ex}%
                                   {\normalfont\large\scshape}}

\newcommand{\R}{\mathbb{R}}
\newcommand{\C}{\mathbb{C}}
\newcommand{\Z}{\mathbb{Z}}
\newcommand{\M}{\mathcal{M}}
\newcommand{\G}{\mathcal{G}}
\newcommand{\F}{\mathscr{F}}
\newcommand{\g}{\Gamma}
\makeatother

\title{ \vspace{-5ex}\bf \large Regular coverings and parallel products of Farey maps}

\author{\normalsize Margaret Stanier\footnote{School of Mathematics and Statistics, The Open University, Milton Keynes, MK7 6AA, United Kingdom \today}}

\date{\vspace{-5ex}}

\begin{document}

\maketitle
\begin{abstract}We examine the structure of Farey maps, which are a class of maps (graph embeddings on surfaces) that have received significant attention recently.  We describe how they are related to each other through regular coverings and parallel products, and use these observations to find their complete spectra, recovering some known results.  We then examine a similar class of maps defined by Hecke groups. 
\end{abstract}



%%%%%%%%%%
\section{Introduction}
%%%%%%%%%%
 Recently, there has been significant research on a class of maps on surfaces known as \emph{Farey maps}, see 
 \cite{IS,SS16,I19,SS,N,KaSi}.    In this paper we use parallel products and regular coverings to describe how they are related to each other, and then calculate their spectra, recovering known results from \cite{N}.  We also obtain a result given in \cite{DLM, Gun} in the context of finite fields: the underlying graphs  of Farey maps of prime level are Ramanujan graphs.   Motivated by this, we study in the same way the maps defined by Hecke groups considered in \cite{IS,KaSi}, which are similar to Farey maps but with faces which are quadrilaterals and hexagons rather than triangles. 

Using the theory developed in \cite{TMOS} and summarised in \cite{Hu,Hu16}, we define an \emph{algebraic orientable regular map} (in what follows we refer to this as a regular map) as a triple
$\M=(G,x,y)$, where $G$ is a finite group, $G=\langle x,y\rangle$, and $(xy)^2=e$, the group identity.  $\M$ is said to be of \emph{type} $(m,n)$, where $m$ and $n$ are, respectively, the orders of $x$ and $y$ in $G$.  A \emph{topological map} is the embedding of a graph $\G$ (the \emph{underlying graph}) on a compact surface $\mathcal{S}$ (the \emph{supporting surface}) such that each component of the complement of $\G$ in $\mathcal{S}$ is homeomorphic to an open disc.  From $\M$ we can construct an oriented topological map with edges defined as the left cosets of $\langle xy \rangle$ in $G$,  vertices as the left cosets of $\langle y\rangle$ in $G$, and faces as the left cosets of $\langle x\rangle$ in $G$; incidence is determined by non-empty intersection. Each member of $G$ is thus identified with one of the directed edges, or \emph{darts}, of this topological representation.  Each vertex is incident to $n$ darts, so $n$ is the vertex valency, and each face is incident to $m$ darts, so $m$ is the face valency.  The vertices and edges determine the underlying graph and the faces the supporting surface, which is oriented in such a way that the $n$ incident darts around each vertex are positioned anticlockwise in the order of increasing powers of $y$.  

In \cite{UT}, the definition of an algebraic map is extended to include a map corresponding to the infinite modular group $\Gamma=\textnormal{PSL}_2(\mathbb{Z})$, which is the quotient by $\{\pm I\}$, where $I$ is the $2\times2$ identity matrix, of the group \[\textnormal{SL}_2(\mathbb{Z})=\left\{\begin{pmatrix}a&b\\c&d\end{pmatrix}:a,b,c,d\in\mathbb{Z}, ad-bc=1\right\}.\]
 We write a member of $\Gamma$ as \[\begin{bmatrix}a&b\\c&d\end{bmatrix}=\left\{\begin{pmatrix}a&b\\c&d\end{pmatrix},\begin{pmatrix}-a&-b\\-c&-d\end{pmatrix}\right\},\quad\text{where}\quad\begin{pmatrix}a&b\\c&d\end{pmatrix}\in\textnormal{SL}_2(\mathbb{Z}).\]
$\g$ is generated by $X=\begin{bmatrix}0&1\\-1&1\end{bmatrix}$, of order $3$, and $Y=\begin{bmatrix}1&1\\0&1\end{bmatrix}$;  $XY=\begin{bmatrix}0&1\\-1&0\end{bmatrix}$ is of order $2$.

\begin{figure}[h]
\begin{center}
\begin{tikzpicture}[scale=0.8]
\draw(0,0)--(14,0);
\draw(7,4)--(7,0);
\draw (1,4)--(1,0);
\draw (13,4)--(13,0);
\draw (1,0) arc [radius=1, start angle=0, end angle=60];
\draw(1,0) arc[radius=3, start angle=180, end angle=0];
\draw(1,0) arc[radius=1, start angle=180, end angle=0];
\draw(1,0) arc[radius=1.5, start angle=180, end angle=0];
\draw(3,0) arc[radius=0.5, start angle=180, end angle=0];
\draw(4,0) arc[radius=1.5, start angle=180, end angle=0];
\draw(4,0) arc[radius=0.5, start angle=180, end angle=0];
\draw(5,0) arc[radius=1, start angle=180, end angle=0];

\draw(7,0) arc[radius=3, start angle=180, end angle=0];
\draw(7,0) arc[radius=1, start angle=180, end angle=0];
\draw(7,0) arc[radius=1.5, start angle=180, end angle=0];
\draw(9,0) arc[radius=0.5, start angle=180, end angle=0];
\draw(10,0) arc[radius=1.5, start angle=180, end angle=0];
\draw(10,0) arc[radius=0.5, start angle=180, end angle=0];
\draw(11,0) arc[radius=1, start angle=180, end angle=0];
\draw (13,0) arc [radius=1, start angle=180, end angle=120];

\node[below] at (9,0){$\frac{1}{3}$};
\node[below] at (10,0){$\frac{1}{2}$};
\node[below] at (11,0){$\frac{2}{3}$};
\node[below] at (13,0){$1$};
\node[below] at (7,0){$0$};
\node[below] at (5,0){$-\frac{1}{3}\phantom{-}$};
\node[below] at (4,0){$-\frac{1}{2}\phantom{-}$};
\node[below] at (3,0){$-\frac{2}{3}\phantom{-}$};
\node[below] at (1,0){$-1\phantom{-}$};

\end{tikzpicture}
\caption{Part of the Farey tessellation}
\label{Fareytessellation}
\end{center}
\end{figure}
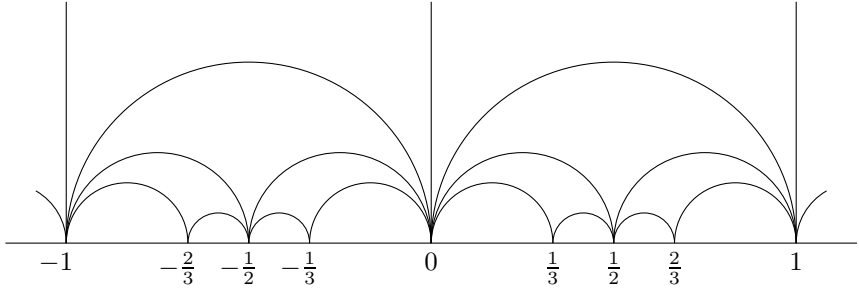

We define an infinite regular algebraic map $\mathscr{F}=(\g,X,Y)$, of type $(3,\infty)$.
In the topological representation of $\mathscr{F}$, any $\gamma\in\Gamma$ is a dart, and the vertices are the cosets \[\gamma\langle Y\rangle=\left\{\begin{bmatrix}a&ar+b\\c&cr+d\end{bmatrix}:r\in\mathbb{Z}\right\}.\]Given $a,c\in \Z$ such that gcd$(a,c)=1$, we can determine a unique vertex of $\F$: we find $b,d\in\Z$ such that $ad-bc=1.$  Then, if $\gamma\in\g$ is the matrix with entries $a,b,c$ and $d$, the vertex corresponding to the ordered pair $(a,c)$ is $\gamma\langle Y\rangle$. Hence we can identify the vertices with the reduced rationals $a/c$, with the usual convention that $1/0=\infty$.  The coset $\gamma\langle XY \rangle$ is an
edge consisting of the two darts $\gamma$ and $\gamma XY$.  We say that $a/c$ is the initial vertex of $\gamma$ and that, as $b/d$ is the initial vertex of $\gamma XY$, it is the final vertex of $\gamma$.  The reduced rationals $a/c$ and $b/d$ are vertices incident to the same edge, or \emph{adjacent} vertices, if and only if $ad-bc=\pm 1$.  $\mathscr{F}$ is often drawn as a tessellation of the upper half plane known as the \emph{Farey tessellation}, shown in Figure \ref{Fareytessellation}.

%%%%%%%%%%%%%%%%%%%%%%%%%%%%%%%%%%%%%%%%%%%
It is shown in \cite[Theorem 1]{UT} that any regular map on an oriented surface with triangular faces is the quotient  of $\F$ by a normal subgroup of finite index of the modular group.  In this sense, $\F$ is the \emph{universal triangular map}.

The {\it  principal congruence subgroup $\g(n)$ of level $n$} is the normal subgroup of $\g$ given by:
 \[\Gamma(n)=\left\{\begin{bmatrix}a&b\\c&d\end{bmatrix}\in \textnormal{PSL}_2(\mathbb{Z}):\begin{bmatrix}a&b\\c&d\end{bmatrix}\equiv \begin{bmatrix}1&0\\0&1\end{bmatrix}\!\!\!\!\!\!\pmod{n}\right \}.\]
%%%%%%%%%%%%%%%%%%%%%%%%%%%%%%%%%

The {\it Farey map of level $n$} is the regular map
$\M_3(n)=(\Gamma/\g(n),X\g(n),Y\g(n))$ of type $(3,n)$.  The group $\g/\g(n)$ is isomorphic to PSL$_2(\Z/n\Z)$.   Its members, the darts of $\M_3(n)$, are, for $\gamma\in\g$,  the cosets $\gamma\g(n)$.  As we wish to compare Farey maps of different levels, we will use the following notation for the members of PSL$_2(\Z/n\Z)$:
\[\begin{bmatrix}a&b\\c&d\end{bmatrix}_n=\left\{\begin{bmatrix}a'&b'\\c'&d'\end{bmatrix}\in \textnormal{PSL}_2(\mathbb{Z}): \begin{bmatrix}a'&b'\\c'&d'\end{bmatrix}\equiv \begin{bmatrix}a&b\\c&d\end{bmatrix}\!\!\!\!\!\pmod{n}\right\}\!.\] 
%%%%%%%%%%%%%%%%%%%%%%%%%%%%%%%%%%%%%%%%%%%%%%%%%%%%%%%%
Recall from \cite{SS} that the subgroup $\g_1(n)$ of $\g$ is given by
\[\g_1(n)=\left\{\begin{bmatrix}a&b\\c&d\end{bmatrix}_n\in\textnormal{ PSL}_2(\Z/n\Z):\quad\begin{bmatrix}a&b\\c&d\end{bmatrix}_n=\begin{bmatrix}1&r\\0&1\end{bmatrix}_n:\quad r=0,1,\ldots,n-1.\right \}\]
Then we note that $\g_1(n)$ is isomorphic to the subgroup of $\Gamma$ generated by $Y\g(n)$, and so the vertices of $\M_3(n)$ are, for $\gamma\in \g$, the cosets $\gamma\g(n)\langle Y\g(n)\rangle=\gamma\langle Y\g(n)\rangle=\gamma\g_1(n)$, or 
 \[\gamma\g_1(n)=\left\{\begin{bmatrix} a&ar+b\\c&cr+d\end{bmatrix}_n:\quad r=0,1,\ldots,n-1.\right\}.\]
Given $a,c\in \Z$ such that gcd$(a,c,n)=1$, we can determine a unique vertex of $\M_3(n)$: we find $b,d\in\Z$ such that $ad-bc\equiv1\!\!\pmod{n}.$  Then, if $\gamma\in\g$ is the matrix with entries $a,b,c$ and $d$, the vertex corresponding to the ordered pair $(a,c)$ is $\gamma\g_1(n)$, which we denote by $[a/c]_n$. 
%%%%%%%%%%%%%%%%%%%%%%%%%%%%%%%%%%%%%%%%%%%%
There is a bijection between these vertices and the \emph{Farey fractions} defined in \cite{SS}, which are the equivalence classes $\{(a',c')\in\Z\times\Z: \textnormal{gcd}(a',c',n)=1,\;(a',c')\equiv \pm(a,c)\!\!\pmod{n}\}$.  Note that $[2/2]_5$ is a vertex of $\M_3(5)$ different to $[1/1]_5$, and $[2/0]_5$ is not the same vertex as $[1/0]_5$. 
 
 An edge of $\M_3(n)$ is the coset $\gamma\g(n)\langle XY \rangle$,
 which consists of the two darts $\gamma\g(n)$ and $\gamma\g(n) XY$.  If $\gamma\in\g$ is the matrix with entries $a,b,c$ and $d$, we say that $[a/c]_n$ is the initial vertex of 
$\gamma\g(n)$ and that $[b/d]_n$  is the final vertex of $\gamma\g(n)$.
  The vertices $[a/c]_n$ and $[b/d]_n$ are adjacent in $\M_3(n)$ if and only if $ad-bc\equiv\pm 1\!\!\!\pmod{n}$.
 
The regular map $\M_1=(G_1,x_1,y_1)$ of type $(m_1,n_1)$ is a \emph{regular covering} of a regular map $\M_2=(G_2,x_2,y_2)$ 
of type $(m_2,n_2)$ if there is a group epimorphism  $\sigma$,  the \emph {covering transformation}, from $G_1$ to $G_2$, with
 $\sigma(x_1)=x_2$ and $\sigma(y_1)=y_2$. If $n_1> n_2$, we have  $\sigma((y_1)^{n_1})=(\sigma(y_1))^{n_1}=y_2^{n_1}$; 
 but, if $e_1$ and $e_2$ are, respectively, the identities of $G_1$ and $G_2$, $\sigma((y_1)^{n_1})=\sigma(e_1)=e_2$.  
Therefore $y_2^{n_1}=e_2$. But $n_2$ is the smallest integer such that $(y_2)^{n_2}=e_2$, so $n_2$ is a divisor of $n_1$. 
 We say that the regular covering is \emph{ramified} at the vertices with \emph{vertex ramification index} $n_1/n_2$.  If $m_1> m_2$, then, similarly,  $m_2$ is a divisor of $m_1$, and the covering is \emph{ramified} at the faces with \emph{face ramification index} $m_1/m_2$. A regular covering of algebraic maps induces a  covering of their supporting surfaces. The \emph{fibre} of a dart $h_2$ of $\M_2$ under a covering transformation $\sigma$ is ${\sigma}^{-1}(h_2)$. The kernel of $\sigma$ in $G_1$, $\textnormal{Ker}(\sigma)$, is the transitive automorphism group of each fibre.  The covering  has $r$ \emph{sheets} if $r=|\textnormal{Ker}(\sigma)|$ (and is a \emph{double covering} if $r=2$). 

\begin{figure}[h]
\begin{center}
\begin{tikzpicture}[scale=0.8]
\path[fill=lightgray](-6.1,0.65) circle [radius=2.3];
\draw[fill=white](-8,-0.5) to [out=340, in=210] (-4.2,-0.5) to [out=100, in=330](-6,2) to [out=200, in =80](-8,-0.5); 
\node[ above left] at  (-4.2,-0.5){$[0/1]_2$};
\node [below] at (-6,1.8){$[1/0]_2$};
\node[above right] at (-8,-0.5){$[1/1]_2$};
\draw [fill](-4.2,-0.5) circle [radius=0.07];
\draw [fill] (-6,2) circle [radius=0.07];
\draw [fill] (-8,-0.5) circle [radius=0.07];

%%%%%%%%%%%%%%%%%%%%

\path[fill=lightgray](-1.5,1)--(0,0)--(1.25,3.25)--(-1.5,1);
\path[fill=lightgray](-1.5,1)--(2.5,1)--(1.25,-2.25)--(-1.5,1);
\path[fill=white](-1.5,1)--(0,0)--(1.25,-2.25)--(-1.5,1);
\path[fill=lightgray](4,0)--(0,0)--(1.25,-2.25)--(4,0);
\path[fill=lightgray](4,0)--(2.5,1)--(1.25,3.25)--(4,0);

\draw (1.25,-2.25)--(0,0);
\draw (1.25,-2.25)--(4,0);
\draw[dashed] (1.25,-2.25)--(2.5,1);
\draw (1.25,-2.25)--(-1.5,1);
\draw [dashed](1.25,3.25)--(2.5,1);
\draw (1.25,3.25)--(-1.5,1);
\draw (1.25,3.25)--(0,0);
\draw (1.25,3.25)--(4,0);
\draw (0,0)--(4,0);
\draw (0,0)--(-1.5,1);
\draw[dashed] (2.5,1)--(-1.5,1);
\draw[dashed] (2.5,1)--(4,0);

\draw [fill] (0,0) circle [radius=0.07];
\draw [fill] (4,0) circle [radius=0.07];
\draw [fill] (2.5,1) circle [radius=0.07];
\draw [fill] (-1.5,1) circle [radius=0.07];
\draw [fill] (1.25,3.25) circle [radius=0.07];
\draw [fill] (1.25,-2.25) circle [radius=0.07];

\node[ left] at (-0.5,0){$[2/1]_4$};
\node[right] at (4,0){$[3/1]_4$};
\node[left] at (-1.5,1){$[1/1]_4$};
\node[above left] at (2.4,1){$[0/1]_4$};
\node[above] at (1.25,3.25){$[1/0]_4$};
\node[below] at (1.25,-2.25){$[1/2]_4$};

\end{tikzpicture}
\caption{The Farey maps $\M_3(2)$ and $\M_3(4)$.}
\label{maps}
\end{center}
\end{figure}
For example, Figure \ref{maps} shows $\M_3(2)$ and $\M_3(4)$.  The underlying surface of each of these maps is a sphere.   $\M_3(4)$ is a four sheeted covering of $\M_3(2)$, ramified at all the vertices with ramification index 2. 
 The covering transformation takes $X\g(4)$ to $X\g(2),$ and $Y\g(4)$ to $Y\g(2)$. 
 So we have $\gamma\g(4)\longrightarrow \gamma\g(2)$ and $[a/c]_4\longrightarrow [a/c]_2$.  For instance
\[\begin{bmatrix}1&2\\0&1\end{bmatrix}_4\longrightarrow \begin{bmatrix}1&0\\0&1\end{bmatrix}_2,\quad \quad[1/0]_4\longrightarrow [1/0]_2,\;\text{ and }\quad[2/1]_4\longrightarrow [0/1]_2.\]

  If ${\G}_1$ and ${\G}_2$ are graphs with dart sets ${\Omega}_1$ and ${\Omega}_2$, respectively, the {\it graph tensor product} $\G_1\times\G_2$ is the graph whose dart set is the cartesian product ${\Omega}_1\times{\Omega}_2$.  If $l$ and $m$ are coprime integers, \cite[Theorem 2]{SW2002} shows that, if $\G_3(l)$ and $\G_3(m)$ are the underlying graphs of $\M_3(l)$ and $\M_3(m)$,  there is another map, the \emph{parallel product } of  the maps $\M_3(l)$ and $\M_3(m)$, whose underlying graph is is $\G_3(l)\times\G_3(m)$.  We will define the parallel product in Section~\ref{3}.  It was introduced in \cite{SW1994}, and used recently in \cite{Hu}.  It is analogous to the join of hypermaps defined in \cite{BrNe}, and to the blend of polytopes used in \cite{MMW}.  
We will prove the following two theorems in Sections 2 and 3.
\begin{theorem}\label{A1}
 For a prime $p$ and a positive integer $k$, $\M_3(p^{k})$ is  a regular covering of $\M_3(p^{k-1})$, ramified at the vertices, with ramification index $p$.  The covering has $4$ sheets if $p^k=4$, and otherwise has $p^3$ sheets.
\end{theorem}
\begin{theorem}\label{A2} Let $m$ be a positive integer. \vspace{-2ex}
\begin{enumerate}[label=\emph{(\roman*)}]
\item If $m$ is odd, $\M_3(2m)$  is the parallel product of the maps $\M_3(2)$ and $\M_3(m)$.
\item If $l$ and $m$ are coprime integers, and neither $l$ nor $m$ is twice an odd integer, then $\M_3(lm)$ is a regular unramified double covering of the parallel product of $\M_3(l)$ and $\M_3(m)$.
\end{enumerate}  
\end{theorem}
A \emph{graph homomorphism} $\pi$ from a graph $\G_1$ to a graph $\G_2$ is a mapping from the vertex set of $\G_1$ to the vertex set of $\G_2$ which preserves adjacency.  The homomorphism $\pi$ is a \emph{graph covering transformation}
if there is bijection between the set of vertices adjacent to a vertex $v$ of $\G_2$ and the set of vertices adjacent to any vertex of $\G_1$ in ${\pi}^{-1}(v)$.  $\G_1$ is then a \emph{graph covering} of $\G_2$, as defined in \cite{GR}.  Note that, if $\G_1$ and $\G_2$ are the underlying graphs of two maps, $\pi$ may not be induced by a map covering transformation.  We will show that, for a prime $p$, $\G_3(p)$ is a graph covering of the complete graph on $p+1$ vertices. 

The eigenvalues of the adjacency matrix of graph $\G$ with $N$ vertices are
 ${\lambda}_1\ge{\lambda}_2\ge\cdots\ge{\lambda}_N$. The {\it spectrum} of $\G$ is the multiset of these eigenvalues with their multiplicities, which we write 
 $\textnormal{sp}(\G)=\{{\lambda}_1^{(m_1)},{\lambda}_2^{(m_2)},\ldots,{\lambda}_i^{(m_i)}\}$ 
if its $i$ distinct eigenvalues are ${\lambda}_1,{\lambda}_2,\ldots,{\lambda}_i$
 with multiplicities $m_1,m_2,\ldots,m_i$ respectively.  If $\G$ is the underlying graph of a map $\M$, the spectrum of $\M$ is defined as $\textnormal{sp}(\M)=\textnormal{sp}(\G)$.  We simplify notation by writing $\textnormal{sp}_3(n)=\textnormal{sp}(\M_3(n))$. 
 We define a product of multisets: if $\textnormal{sp}(\M_1)=\{{\lambda}_1^{(m_1)},{\lambda}_2^{(m_2)},\ldots,{\lambda}_i^{(m_i)}\}$, and $\textnormal{sp}(\M_2)=\{{\mu}_1^{(l_1)},{\mu}_2^{(l_2)},\ldots,{\mu}_j^{(l_j)}\}$, then
$\textnormal{sp}(\M_1)\textnormal{sp}(\M_2)=\{{\lambda}_r{\mu}_s^{(m_rl_s)}:\;r=1,\ldots,i;\;s=1,\ldots,j\}.$  Also, for   $k\in\mathbb{Z}$, we write $k\textnormal{sp}(\M_1)=\{k{\lambda}_1^{(m_1)},k{\lambda}_2^{(m_2)},\ldots,k{\lambda}_i^{(m_i)}\}.$ 

  In \cite{DFMRS,DFPS,DFS}  the spectra of graphs are found using coverings ramified at the face centres of maps.
 We use Theorems~\ref{A1} and \ref{A2}, and regular coverings ramified at map vertices, to find, iteratively, the spectrum of the Farey map $\M_3(n)$ from the prime decomposition of $n$.  
We have
 \[\textnormal{sp}_3(2)=\{-1^{(2)},2\},\quad\textnormal{sp}_3(3)=\{-1^{(3)},3\}\;\text{ and }\;\textnormal{sp}_3(4)=\{-2^{(2)},0^{(3)},4\}.\]  For higher values of $n$, we prove the following two theorems in Sections 6 and 7.  We need the formula for the  number of vertices of $\M_3(n)$, which is, from \cite{IS}, 
\begin{eqnarray}|V_3(n)|&=&\frac{n^2}{2}\prod_{p|n}\left(1-\frac{1}{p^2}\right),
\end{eqnarray}
where the product is over all prime divisors of $n$.
\begin{theorem}\quad\label{B1} Let $p$ be a prime, and $k$ a positive integer. \vspace{-1ex}
\begin{enumerate}[label=\emph{(\roman*)}]
\item If $p>3$, then
\vspace{-2 ex}\[\textnormal{sp}_3(p)=\{-\sqrt{p}^{(m)}\!\!,-1^{(p)}\!\!,\sqrt{p}^{(m)}\!\!,p\},\;\]$\text{where }\;m=\frac{1}{4}(p-3)(p+1).$
\item If $p^k>4$, then 
\vspace{-2 ex}\[\textnormal{sp}_3(p^{k})=p\ \!\textnormal{sp}_3(p^{k-1})\cup \{-\sqrt{p^{k}}^{(\frac{1}{2}pc)}\!\!,  0^{(c)}\!\!, \sqrt{p^{k}}^{(\frac{1}{2}pc)}\},\]$\text{  where }c=(p-1)|V_3(p^{k-1})|.$
\end{enumerate}
\end{theorem}
\noindent For example, because $\textnormal{sp}_3(7)=\{-\sqrt{7}^{(8)}\!\!,-1^{(7)}\!\!,\sqrt{7}^{(8)}\!\!,7\}$, it follows that
\[\textnormal{sp}_3(49)=\{-7\sqrt{7}^{(8)}\!\!,-7^{(511)}\!,0^{(144)}\!,7^{(504)}\!,7\sqrt{7}^{(8)}\!\!,49\}.\]
%%%%%%%%%%%%%%%%%%%%%%%%%%%%%%%%%%%%%%%
\begin{theorem}\label{B2}\; Let $m$ be a positive integer. \vspace{-1ex}
\begin{enumerate}[label=\emph{(\roman*)}]
\item If $m$ is odd, then $\textnormal{sp}_3(2m)=\textnormal{sp}_3(2)\textnormal{sp}_3(m)$.
\item If $l$ and $m$ are coprime integers, and neither $l$ nor $m$ is twice an odd integer, then
\vspace{-1ex} \[\textnormal{sp}_3(lm)=\textnormal{sp}_3(l)\textnormal{sp}_3(m)\cup\{-\sqrt{lm}^{(N/4)},\sqrt{lm}^{(N/4)}\}, \]$\text{ where }N=|V_3(lm)|.$
\end{enumerate}
\end{theorem}
For example,\vspace{-2ex}
\begin{eqnarray*}
\textnormal{sp}_3(14)&=&\{-1^{(2)}\!,2\}\{-\sqrt{7}^{(8)}\!\!,-1^{(7)}\!, \sqrt{7}^{(8)}\!\!,7\}\\
&=&
\{-7^{(2)}\!,-2\sqrt{7}^{(8)}\!\!,-2^{(7)}\!,-\sqrt{7}^{(16)}\!\!,1^{(14)}\!,\sqrt{7}^{(16)}\!\!,2\sqrt{7}^{(8)}\!\!,14\}.\\
\textnormal{sp}_3(28)&=&\{-2^{(2)}\!,0^{(3)}\!,4\}\{-\sqrt{7}^{(8)}\!\!,-1^{(7)}\!, \sqrt{7}^{(8)}\!\!,7\}\cup\{-2\sqrt{7}^{(72)}\!\!,2\sqrt{7}^{(72)}\}\\
&=&
\{-14^{(2)}\!,-4\sqrt{7}^{(8)}\!\!,-2\sqrt{7}^{(88)}\!\!,-4^{(7)}\!,0^{(72)}\!,2^{(14)}\!,2\sqrt{7}^{(88)}\!\!,4\sqrt{7}^{(8)}\!\!,28\}.
\end{eqnarray*}
We also show that the underlying graphs of certain maps defined by Hecke groups are graph double coverings of the underlying graphs of the Farey map $\M_3(n)$.  We use this result to find their spectra.
%%%%%%%%%%%%%%%%%%%%%%%%%%%%%%%%%%%%%%%%%%%%%%%%
\section{Farey maps and regular coverings:  the proof of Theorem~\ref{A1}.}\label{2}
We prove Theorem~\ref{A1} by giving a more general result, which we will also use later.
\begin{lemma}\label{d}If $n=dm$, the Farey map $\M_3(n)$ is a regular map covering of $\M_3(m)$ of order $d^3$ if $d\neq 2$, or $4$ if $d=2$, which is ramified at the vertices with ramification index $d.$ The covering transformation takes the dart $\gamma\g(n)$ to the dart $\gamma\g(m)$, and the vertex $[a/c]_{n}$ to the vertex $[a/c]_{m}$.
\end{lemma}
\begin{proof}
We have $\M_3(n)=(\g/\g(n),X\g(n),Y\g(n))$  and $\M_3(m)=(\g/\g(m),X\g(m),Y\g(m))$. We define the mapping
\begin{eqnarray*}\sigma:\g/\g(n)&\longrightarrow&\g/\g(m)\\
X\g(n)&\longmapsto& X\g(m);\\Y\g(n)&\longmapsto &Y \g(m),
\end{eqnarray*}
This is is an epimorphism and so a covering transformation, since, as $m$ is a divisor of $n$, $\g(n)$ is a subset of $\g(m)$ .  As
 \[\sigma\left(\begin{bmatrix}1+km&rm\\sm&1-km\end{bmatrix}_{n}\right)=\begin{bmatrix}1&0\\0&1\end{bmatrix}_{m}\;\text{ for }k,r,s=0,1,\ldots,d-1,\]
the kernel of $\sigma$ has order $d^3$ if $d\neq 2$, or 4 if $d=2$.

Since $\g=\langle X,Y \rangle,$ for all $\gamma\in\g,\quad\sigma(\gamma\g(n))=\gamma \g(m)$.  The vertices of $\M_3(n)$ are the cosets $\gamma\g_1(n) =[a/c]_{n}$. So 
$\sigma(\gamma\g_1(n))=\gamma\g_1(m)=[a/c]_{m}.$
\end{proof}%%%%%%%%%%%%%%%%%%%%%%%%%%%%%%%%%%%%%%%
Putting $n=p^{k-1}, d=p$ in Lemma \ref{d} then proves  Theorem \ref{A1}.

\section{Farey maps and parallel products: the proof of Theorem~\ref{A2}.}\label{3}%%%%%%%%%%%%%%%%%%%%%%%%%%%%%%%%%%%%%%%%%%%%%%%%%%%%%%%%%%%%%%%%%%%%%%%%%%%%
 In order to decompose $\M_3(n)$ for a composite $n$, we use a product of maps which is consistent with the tensor product of their underlying graphs, following the approach in
 \cite{SW1994}. 
\begin{definition}\label{pp} Let $G_1$ be a finite group generated by $x_1$ and $y_1$, and let $G_2$ a finite group generated by $x_2$ and $y_2$.  We define the {\it parallel product}  of $(G_1,x_1,y_1)$ and $(G_2,x_2,y_2)$ as the group $(E,(x_1,x_2),(y_1,y_2))$, where $E$ is the subgroup of $G_1\times G_2$ generated by $(x_1,x_2)$ and $(y_1,y_2)$.  
\end{definition}
\noindent If $x_1y_1$ and $x_2y_2$ are both of order $2$, so is $(x_1y_1,x_2y_2)$.  Then the triple $(E,(x_1,x_2),(y_1,y_2))$ defines a regular map, which is the parallel product of the maps $\M_1=(G_1,x_1,y_1)$ and $\M_2=(G_2,x_2,y_2)$.  If $G_1=\g/H$ and $G_2=\g/K$, where $H$ and $K$ are both normal subgroups of $\g=\langle X,Y \rangle$, then from \cite[Lemma 3(vii)] {Hu} and \cite[Theorem 3.1]{BrNe},  the parallel product is isomorphic to the map $(\g/H\cap K,XH\cap K,YH\cap K)$.

\noindent  We now prove Theorem \ref{A2}.  It gives an iterative method for decomposing the Farey map $\M_3(n)$ into a series of regular coverings and parallel products given the prime decomposition of $n$.

\begin{proof}[Proof of Theorem \ref{A2}]

Let $l$ and $m$ be positive coprime integers. Then $a\equiv b\pmod{lm}$ if and only if  $a\equiv b\pmod{l}$ and $a\equiv b\pmod{m}$, so that $\g(lm)\subset\g(l)\cap\g(m)$.  So we can define the group epimorphism
\begin{eqnarray*}\sigma:\g/\g(lm)&\longrightarrow&\g/(\g(l)\cap\g(m))\\
\gamma\g(lm)&\longmapsto&\gamma(\g(l)\cap\g(m)).
\end{eqnarray*}
 Then $\sigma$ is a regular map covering transformation.  It maps each member of $\g/\g(lm)$, a coset $\gamma\g(lm)$, onto  the coset $\gamma\g(l)\cap\g(m)$ in which it is contained, which is a member of $\g/\g(l))\cap\g(m)$. In particular the kernel of $\sigma$ consists of those members of $\g/\g(lm)$ which are mapped onto the identity of $\g/\g(l)\cap\g(m)$.  This is the coset of $\g(l)\cap\g(m)$ in $\g$ containing the identity of $\g$, which is $\g(l)\cap\g(m)$.  If a matrix with entries $a,b,c$ and $d$ is a member of $(\g/\g(lm))\cap(\g(l)\cap\g(m))$, then $a\equiv d\equiv\pm1\pmod{l},b\equiv c\equiv0\pmod{l}$, and  $a\equiv d\equiv\pm1\pmod{m},b\equiv c\equiv0\pmod{m}$.   Then $b\equiv c\equiv0\pmod{lm}$ and either $a\equiv d\equiv\pm1\pmod{lm}$ or $a\equiv d\equiv\pm u\pmod{lm}$ for any $u\in U$, where $U=\{u\in\mathbb{Z}: u\equiv 1\pmod{r}\text{ and }u\equiv -1\pmod{s}\}$.  If either $l=2$ or $m=2$ (or both), then $U=\emptyset$,  so that $\g(2)\cap\g(m)=\g(2m)$, and $\M_3(2s)$ is the parallel product of $\M_3(2)$ and $\M_3(m)$.  If neither $l$ nor $m$ is equal to 2,  as gcd$(l,m)=1$, there is exactly one $u$, modulo $lm$, such that $u\equiv 1\pmod{l}$ and $u\equiv -1\pmod{m}$. So, as ${\sigma}^{-1}(\g(l)\cap\g(m))=\{\g(lm),u\g(lm)\}$ is of order 2, $\M_3(lm)$ is a double covering of the parallel product of $\M_3(l)$ and $\M_3(m)$, unramified at the vertices as both maps have vertex valency $lm$.
\end{proof}
%%%%%%%%%%%%%%%%%%%%%%%%%%%%%%%%%%%%%%%%%%
\section{Evaluation of spectra using regular coverings.}\label{gc}%Section 3%
%%%%%%%%%%%%%%%%%%%%%%%%%%%%%%%%%%%%%%
In this section we give some general results for regular graphs and maps.  Recall that, if $\M_1$ is a regular covering of $\M_2$ with covering transformation $\sigma$, and $h$ is a dart of $\M_2$, then $\sigma^{-1}(h)$ is the fibre of $h$.  If $v$ is a vertex of $\M_2$, we define the fibre of $v$ as  $\sigma^{-1}(v)$.
\begin{lemma}\label{covervalency}%3.1%
Let $\M_1=(G_1,x_1,y_1)$ and $\M_2=(G_2,x_2,y_2)$ be regular maps with vertex valencies $n$ and $m$ respectively, and let ${\mathcal{M}}_1$ be a regular covering of ${\mathcal{M}}_2$.  Then if $v$ and $v'$ are adjacent vertices of ${\mathcal{M}}_2$, each of the vertices in the fibre of $v$ is adjacent to exactly $d=n/m$ of the vertices in the fibre of $v'$.
\end{lemma}
\begin{proof} Let $\sigma$ be the covering transformation, and let $r$ be the number of sheets of the covering, so that $r=|\textnormal{Ker}(\sigma)|$. Let $h$ be a dart of $\M_2$ with initial vertex $v$ and final vertex $v'$. Let $|{\sigma}^{-1}(v)|=s,$  and let $w$ be a vertex of $\M_1$ in ${\sigma}^{-1}(v)$. 
Then the orbit of $w$ under the action of $\textnormal{Ker}(\sigma)$ is ${\sigma}^{-1}(v)$. 
 Let $d$ be the number of vertices in ${\sigma}^{-1}(v')$ adjacent to $w$. 
 Then the stabiliser of $w$ in $\textnormal{Ker}(\sigma)$ is the set of darts in ${\sigma}^{-1}(h)$ incident to $w$.  By the orbit stabiliser theorem $r=sd$.
  Now let $D$ be the set of all darts in the fibres of all of the $m$ darts incident to $v$.  We have $|D|=rm=sn$, so $dm=n$.
\end{proof}
\begin{lemma}\label{graphcovering}%3.2%
Let a graph $\G_1$ be a regular covering of a graph $\G_2$ with graph covering transformation $\phi$.  Let $v$ and $v'$ be two adjacent vertices of $\G_2$.  Then any vertex in ${\phi}^{-1}(v)$ is adjacent to exactly one of the vertices in ${\phi}^{-1}(v')$.
\end{lemma}
\begin{proof}
Let $w$ be a member of ${\phi}^{-1}(v)$.  There is a bijection between $v$ and its adjacent vertices and $w$ and its adjacent vertices which takes $v'$ to just one vertex adjacent to $w$.
\end{proof}
The \emph{adjacency matrix} of a graph $\mathcal{G}$ on $N$ vertices, $v_1,v_2,\ldots,v_N$, is an $N\times N$ square matrix whose rows and columns are indexed by the graph vertices.  The entry in the $i^{th}$ row and the $j^{th}$ column is $e_{ij}$. If $v_i$ and $v_j$ are adjacent, $e_{ij}=1$, otherwise $e_{ij}=0$.  An adjacency matrix is symmetrical for an undirected graph, and all its diagonal entries are 0 if the graph contains no loops. By definition,
 $\lambda\in\R$ is an eigenvalue of an $N\times N$ symmetric matrix $C$  with eigenvector $x$ if $Cx=\lambda x$.  If $A$ is an adjacency matrix, an eigenvector $x$ has $N$ components $x_v$ for each of $v=v_1,v_2,\ldots,v_N$, and $\lambda$ is an eigenvalue of $A$  with eigenvector $x$ if and only if, for all $v$, \[{\lambda}x_v=\sum_{u\sim v}x_u,\] where the sum is over all the vertices of the graph adjacent to $v$.
%%%%%%%%%%%%%%%%%%%%%%%%%%%%%%%%%%%%%%%%%%%%%%%%%%%%%%%%%%%%%%%%%%%%%%%%%%% 
\begin{lemma}\label{covereig}%3.3%
Let ${\mathcal{G}}_1$ and ${\mathcal{G}}_2$  be two regular graphs with vertex valencies $n$ and $m$ respectively, where $n=dm$, and suppose either that they are the underlying graphs of two regular maps $\M_1$ and $\M_2$, and that  $\M_1$ is a regular map covering of  $\M_2$, or, if $d=1$, that ${\mathcal{G}}_1$ is a graph covering of ${\mathcal{G}}_2$. Let $N$ be the number of vertices of $\mathcal{G}_2$.  Then
\begin{enumerate}[label=\emph{(\roman*)}]
\item $d\ \!\textnormal{sp}(\G_2)\subset\textnormal{sp}(\G_1)$
\item if $d>2$, $d\ \!\textnormal{sp}(\G_2)\cup\{0^{(\gamma)}\}\subset\textnormal{sp}(\G_1)$; $\gamma\ge \begin{cases}N(d-1)\text{ if }d\text{ is odd,}\\  N(d-2)\text{ if }d\text{ is even. }\end{cases}$
\end{enumerate} 
\end{lemma}
\begin{proof}
Let $\pi$ be the regular map covering transformation or graph covering transformation of $\G_2$ by $\G_1$.  Let $w$ be a vertex of $\G_1$ and $v$ a vertex of $\G_2$ with $\pi(w)=v$.  Then, from Lemmas~\ref{covervalency} and \ref{graphcovering}, if $v'$ is adjacent to $v$ in ${\mathcal{G}}_2$, $w$ is adjacent to $d$ of the vertices in ${\pi}^{-1}(v')$.   Let ${\lambda}$ be an eigenvalue of ${\mathcal{G}}_2$ with eigenvector $x$.
\begin{enumerate}[label=(\roman*)]
\item  Define the vector $y$ by  $y_w=x_{\pi(w)}$ for all vertices $w$ of ${\mathcal{G}}_1$.  Then 
 \[\sum_{w'\sim w}y_{w'}=d\sum_{v'\sim v}x_{v'}=d{\lambda}x_v=d\lambda y_w.\] 
So $d{\lambda}$ is an eigenvalue of ${\mathcal{G}}_1$ with eigenvector $y$.  Since a different vector $y$ corresponds to each vector $x$ in the eigenspace of $\lambda$, the geometric multiplicity of $d\lambda$ in ${\mathcal{G}}_1$ is at least that of ${\lambda}$ in ${\mathcal{G}}_2$.
\item Let $v$ be a vertex of  ${\mathcal{G}}_2$, and let $w$ be a vertex in the fibre of $v$.  Let the $m$ neighbours of $v$ be $v'_i$ for $i=1,\ldots,m.$  For each $w\in {\pi}^{-1}(v)$, label the $d$ vertices in the fibre of each $v'_i$ adjacent to $w$ as $w'_{i1},\ldots,w'_{id}$. Let $K$ be an integer such that
 $1\le K<d$ for odd $d$, or $1\le K<d-1$ for even $d$.  Then, for any pair $(v,K)$, we define a vector  $z$ by its components $z_u$ 
(here the component $z_u$ corresponds to the vertex $u$ of  ${\mathcal{G}}_1$):
\[z_u=\begin{cases}\phantom{-}1\text{ if } u=w'_{iK}\quad\;\ \ \text{ for some }w\in {\pi}^{-1}(v)\text{ and some }v_i\text{ adjacent to }v,\\-1\text{ if } u=w'_{i(K+1)}\text{ for some }w\in {\pi}^{-1}(v)\text{ and some }v_i\text{ adjacent to }v,\\\phantom{-}0\text{ otherwise}.\end{cases}\] 
Let $A$ be the adjacency matrix of ${\mathcal{G}_1}$.  Then, for each pair $(v,K)$, the $r$'th component of the vector $Az$ is
\[\sum_ue_{ur}z_u=\sum_{w\in F}\left(\sum_{i=1}^me_{w'_{iK},w}-\sum_{i=1}^me_{w'_{i(K+1)},w}\right)=\sum_{w\in F}\sum_{i=1}^m(1-1)=0z_r,\] 
where $F={\pi}^{-1}(v)$.
So $0$ is an eigenvalue of ${\mathcal{G}}_1$ with eigenvector $z$.
The vectors $z$ together with $y$ form a linearly independent set.  So, as $\mathcal{G}_2$ has $N$  vertices $v$, and $K$ takes $d-1$ or $d-2$ values, 0 is an eigenvalue of ${\mathcal{G}}_1$ with (possibly additional) geometric multiplicity at least $N(d-1)$ for odd $d$, or $N(d-2)$ for even $d$.
\end{enumerate}\vspace{-5ex}
\end{proof}
%%%%%%%%%%%%%%%%%%%%%%%%%%%

We will also use the following result.  
\begin{lemma}\label{blocks}
Let $A$ be a symmetric $N\times N$ matrix, and suppose that, for some integer $n$, \begin{eqnarray*}\label{pmatrix}
A^2-nI&=&\begin{pmatrix}C&C&\cdots&C\\C&C&\cdots&C\\\vdots&\vdots&&\vdots\\C&C&\cdots&C\end{pmatrix}, 
\end{eqnarray*} where $C$ is an $r\times r$ symmetric
 matrix. Then either $\sqrt{n}$ or $-\sqrt{n}$ or both are eigenvalues of $A$ with total geometric multiplicity greater than or equal to $N-\textnormal{rank}(C)$.
\end{lemma}
\begin{proof}
The rows and columns of $A^2-nI$ are not linearly independent, so its determinant is 0 and $n$ is an eigenvalue of $A^2$.  Since $|A^2-nI|=|A-\sqrt{n}I||A+\sqrt{n}|$, $\sqrt{n}$ or $-\sqrt{n}$ or both are eigenvalues of $A$ with total multiplicity equal to the multiplicity of $n$ as an eigenvalue of $A^2$.
 The geometric multiplicity of an eigenvalue $\lambda$ of a matrix $M$ is ${\gamma} (\lambda)$.  It is the dimension of the space generated by its eigenvectors, and equal to the nullity of  $M-\lambda I$. So
\begin{eqnarray*}\label{eigen}{\gamma} (\lambda)&=&M-\text{rank}(M-\lambda I).\end{eqnarray*}
The result follows as the rank of $A^2-nI$ is equal to the rank of $C$.
\end{proof}
%%%%%%%%%%%%%%%%%%%%%%%%%%%%%%%%%%
\section{Farey maps and complete graphs}
We show that the underlying graphs of Farey maps of prime level are graph coverings of complete graphs.  The spectra of complete graphs are known, so this will enable us to find the spectra of the Farey maps. We first collect some necessary information about the vertices of $\M_3(n)$.

The {\it poles} of the map $\M_3(n)$ are defined in \cite{SS} as the vertices $[a/0]_n$, with $a$ coprime to $n$. Let $v$ be any vertex of $\M_3(n)$ and let $M$ be an automorphism which takes the vertex $[1,0]_n$ to $v$.  Then we define the {\it copoles} of $v$ as $M[a,0]_n$ for $a=1,\ldots,h$.  If $d$ is such that $ad\equiv \pm 1\!\!\pmod{n}$ and $1\le d \le n/2$, then the vertices $[b/d]_n$ for $b=0,\ldots,n-1$ are adjacent to $[a/0]_n$.  The \emph{star} a vertex consists of that vertex and all vertices adjacent to it.
\cite[Theorem 7]{SS} shows that the stars of the $h=\frac{1}{2}(p-1)$ poles of $\M_3(p)$ for a prime $p$
 each contain $p+1$ vertices, are disjoint, and, together,  include all of the $\frac{1}{2}(p^2-1)=h(p+1)$ vertices of $\M_3(p)$. 

\begin{theorem}\label{Kp+1} The underlying graph of the Farey map $\M_3(p)$, for an odd prime $p$, is a graph covering of order $\frac{1}{2}(p-1)$ of the complete graph $K_{p+1}$ on $p+1$ vertices. 
\end{theorem}
\begin{proof} We label the $p+1$ vertices of $K_{p+1}$ as $0,1,\ldots,p$. We define a transformation $\phi$, which takes each vertex of $\M_3(p)$ together with all its copoles to the same vertex of $K_{p+1}$.
\begin{eqnarray*}
\phi: \text{ vertices of }\M_3(p)&\longrightarrow& \text{ vertices of }\ K_{p+1}\\
\quad[a/0]_p&\longmapsto&p\\
\text{for }b\neq 0,\quad[a/b]_p&\longmapsto&ab^{-1}.
\end{eqnarray*}
There is a bijection $\tau$ between $[1/0]_p$ and its adjacent vertices, and the vertices of $K_{p+1}$:
\begin{eqnarray*}
\tau: \text{ vertices of the star of}\ [1/0]_p&\longrightarrow& \text{ vertices of }\ K_{p+1}\\
\quad[1/0]_p&\longmapsto&p\\
\text{for }b=0,\ldots,p-1,\quad[b/1]_p&\longmapsto&b. 
\end{eqnarray*}
Let $v$ be any vertex of $\M_3(p)$, and let $M$ be an automorphism of $\M_3(p)$ which takes the vertex $[1/0]_p$ to $v$. Then the transformation $\tau M$ is a bijection between the star of $v$ and the set of vertices of $K_{p+1}$;  therefore $\phi$ is a graph covering.  It has $h=\frac{1}{2}(p-1)$ sheets as $p$ has $h$ pre-images.
\end{proof} 
%%%%%%%%%%%%%%%%%%%%%%%%%%%%%%%%%%%%%%%%%%%%%%%%%%%
\section{The spectrum of $\M(p^k)$ for a prime $p$: proof of Theorem \ref{B1}.}
%%%%%%%%%%%%%%%%%%%%%%%%%%%%%%%%%%%%%%%%%%%%%%%%%
We denote the adjacency matrix of $\M_3(n)$ as  $A(n)$, and index its rows and columns according to the labelling we choose for the vertices of $\M_3(n)$. In particular we will label the vertex $[1/0]_n$ as $v_0$, so that it corresponds to the first row and column of $A(n)$.  The entry of $A(n)^2$ corresponding to the vertices $v_i$ and $v_j$ is equal to the number of walks of length 2 connecting them. (See, for instance, \cite[Lemma 8.1.2]{GR}).  We find this by extending the results given in \cite[Theorems 11-15]{SS} for the distance between two Farey map vertices. 
\begin{lemma}\label{a2entries} Let $b,d$ be such that $\textnormal{gcd}(b,d,n)=1$, so that $v_l=[b/d]_n$ is a vertex of $\M_3(n)$, and let $\textnormal{gcd}(d,n)=r$. Then the number of walks of length $2$ between $v_0$ and $v_l$ is $2$ if $r=1$, $0$ if $r$ is not a divisor of either of $b\pm1$, $r$ if $r$ is a divisor of one of $b\pm1$, and $4$ if $r=2$ and $b$ is odd.
\end{lemma}
\begin{proof}
  Since all vertices adjacent to
 $[1/0]_n$ are $[x/1]_n$ for $x\in{\mathbb{Z}}$, the vertex $v_{l}$ has a walk of length 2 to $[1/0]_n$
 if and only if there is a vertex $[x/1]_n$ adjacent to $[b/d]_n$, that is if and only if there is an integer $x$ such that
\[
xd \equiv b\pm1\!\!\!\pmod{n}.\]
Then our result is the total number of solutions modulo $n$ to these 2 congruences.  We prove it by recalling, for instance from \cite[Theorems 5.12-5.14]{A}, that each congruence has one solution modulo $n$ if and only if gcd$(d,n)=1$, no solutions if gcd$(d,n)=r$ and $r$ does not divide either of $d\pm 1$, and $r$ solutions if $r\neq 2$ and $r$ divides either $b+1$ or $b-1$.  If $r=2$ and $b$ is odd, $2$ divides both $b+1$ and $b-1$, giving $4$ walks of length $2$. 
 \end{proof}
Given any two vertices $v_i=[a/c]_n$ and $v_j=[b/d]_n$ of $\M_3(n)$ we will use this lemma to find the $ij$ entry of $A(n)^2$ in the following way:  we define $\Delta_n(ij)=ad-bc$, and let $\lambda_n(i),\mu_n(i)\in\Z$ be any two integers such that $\lambda_n(i)(a)+\mu_n(i)(c)+\nu n=1$ for some $\nu\in\Z$. Then there is an automorphism $M_n(i)$ which takes $v_i$ to $v_0$, and $v_j$ to $v_l=[\beta_n(ij)/\Delta_n(ij)]_n$, where  $\beta_n(ij)=\lambda_n(i)b+\mu_n(i)d$. The $ij$ entry of $A(n)^2$  is then given by  Lemma~\ref{a2entries}.  
\begin{lemma}If $\Delta_n(ij)\equiv0\!\!\pmod{n}$, the $ij$ entry of the matrix $A(n)^2-nI$ is $0$.\end{lemma}\begin{proof}
If $\Delta_n(ij)=\nu n$ for some $\nu\in\Z$, $M_n(i)$ takes $v_i$ and $v_j$ to $[1/0]_n$ and $[\beta_n(ij)/\nu n]_n$.
 Then $r=n$ divides $\beta_n(ij)\pm1$ if and only if $\beta_n(ij)=\pm1+\kappa n$ for some $\kappa\in\Z$,
 in which case $[\beta_n(ij)/\nu n]_n=v_0$; so $v_i=v_j$, and the $ii$ entry of $A(n)^2$ is $n$, as expected, since there are
 $n$ paths to and from $v_i$ along each of the $n$ edges incident to $v_i$. 
 If $\beta_n(ij)\not\equiv\pm 1\!\!\pmod{n}$, $v_i\neq v_j$.  In this case the $ij$ entry of $A(n)^2$ is 0, as  
 \cite[Theorem 14]{SS} shows that the shortest path between these vertices is of length 3.
\end{proof}
The underlying graphs of ${\mathcal{M}}_3(2)$ and ${\mathcal{M}}_3(3)$ are the complete graphs on 3 and 4 vertices respectively.  The spectrum of the complete graph on $k+1$ vertices is $\{-1^{(k)},k\}$, therefore sp$_3(2)=\{-1^{(2)},2\}$ and sp$_3(3)=\{-1^{(3)},3\}$.   We can now prove 
Theorem \ref{B1}(i).

\begin{proof}[Proof of Theorem \ref{B1}(i)]

 All congruences in this proof are modulo $p$.  As the spectrum of $K_{p+1}$ is $\{-1^{(p)},p\}$, from Lemmas \ref{covereig} and \ref{Kp+1}, $\{-1^{(p)},p\}\subset\textnormal{sp}(p)$. 
  If $\Delta_p(ij)\not\equiv0$, gcd$(\Delta_p(ij),p)=1 $ and the $ij$ entry of $A(p)^2$ is 2. 
 We order the vertices of $\M_3(p)$ so that, if $v_i=[a/c]_p$, $i=a(p+1)$ when $c\equiv0$, and $i=c^{-1}p+ac^{-1}$ if $c\not\equiv0$. 
 This ensures that the $p+1$ vertices of each of the $h=\frac{1}{2}(p-1)$ stars are together, with the poles $p+1$ positions apart. 
Then $\Delta_p(ij)\equiv0$ if, putting $v_j=[b/d]_n$, $ad-bc\equiv0$.  
If $c\equiv0$, as $a\not\equiv0, d\equiv0$, so $i=a(p+1)$ and $j=b(p+1)$. 
 If $c\not\equiv0, i=c^{-1}p+ac^{-1}, j=d^{-1}p+bd^{-1}$, and $ac^{-1}- bd^{-1}\equiv0$,
 so the $ij$ entry of $A(p)^2$ is $0$  if and only if $i$ and $j$ are a multiple of $p+1$ positions apart. Therefore
 \[A(p)^2-pI=\begin{pmatrix}C&C&\cdots&C\\C&C&\cdots&C\\\vdots&\vdots&&\vdots\\C&C&\cdots&C\end{pmatrix}\quad\text{where}\quad
 C=\begin{pmatrix}0&2&\cdots&2\\2&0&\cdots&2\\ \vdots&\vdots&&\vdots\\2&2&\cdots&0\end{pmatrix}.\]  
The matrix $A(p)^2-pI$ is arranged as $h\times h$ copies of the $(p+1)\times(p+1)$ matrix $C$.  The  $ij$ element of $C$ is equal to $2$ if $i\neq j$, or to $0$ if $i=j$.  Then, from Lemma \ref{blocks}, $\sqrt{p}$ or $-\sqrt{p}$ or both are eigenvalues of $\mathcal{M}_3(p)$ with $\gamma(-\sqrt{p})+\gamma(\sqrt{p})\ge\frac{1}{2}(p+1)(p-3)$. 
 The result follows as the total multiplicity of the eigenvalues of $\mathcal{M}_3(p)$ is the number of its vertices, or $\frac{1}{2}(p+1)(p-1)$.
\end{proof}

To prove Theorem \ref{B1}(ii) we need the following lemmas.  We define $q=p^{k-1}$.
\begin{lemma}\label{delta}Let $\Delta$ and $a$ be integers such that  $\Delta+aq\not\equiv0\!\!\pmod{pq}$.  Then we have \[\textnormal{gcd}(\Delta,q)=\textnormal{gcd}(\Delta+aq,pq).\]
\end{lemma}
\begin{proof}Let $r=\textnormal{gcd}(\Delta,q)$  and $s=\textnormal{gcd}(\Delta+aq,pq)$.  It is straightforward to see that $r\leq s$.   Then, as $\Delta+aq\not\equiv0\!\!\pmod{pq}$, $s=p^l$, where $0\le l\leq k-1$.  Hence $s$ divides $q$, and so, since $s$ divides $\Delta +aq$,  $s$ divides $\Delta$.  Consequently $s\leq r$.  Therefore $s=r, $ as required.
\end{proof}
\begin{lemma}\label{sqrt}
If $p^k>4,$ and the number of vertices of ${\mathcal{M}}_3(p^{k-1})$ is $N$, then
 \[\{(-\sqrt{p^{k}})^{(m)},(\sqrt{p^{k}})^{(m)}\}\subset\textnormal{sp}(p^{k}), \;\text{ where }m=\frac{1}{2}p(p-1)N.\]
\end{lemma}
\begin{proof}
 From Lemma \ref{d}, there is a regular covering of ${\mathcal{M}}_3(q)$ by ${\mathcal{M}}_3(pq)$. Let $v_i=[a/c]_q$ and $v_j=[b/d]_q$ be vertices of $\M_3(q)$.   Let $w_f$ and $w_g$ be vertices of $\M(pq)$ such that $w_f$ is in the fibre of $v_i$,
 so that $w_f=[a+sp/c+tp]_{pq}$ for some $s,t\in\Z$, and  $w_g$ is in the fibre of $v_j$.

We compare the entries of $A(pq)^2$ and $A(q)^2$.  We can check that  $\beta_{pq}(fg)=\beta_{q}(ij)+\rho q,$ and $ \Delta_{pq}(fg)=\Delta_{q}(ij)+\tau q$ for some $\rho,\tau\in\Z$. We first assume that $\Delta_{pq}(fg)\not\equiv0\!\!\pmod{pq}$.  From Lemma \ref{delta}, $\textnormal{gcd}(\Delta_{pq}(fg),pq)=\textnormal{gcd}(\Delta_q(ij),q)=r$. So, from Lemma \ref{a2entries}, the $fg$ entry of $A(pq)^2$ is equal to the $ij$ entry of $A(q)^2$.  We order the rows and columns of the $p^2N\times p^2N$ matrix $A(pq)^2$ so that $w_f=[a+sp/c+tp]_{pq}$ is in position $f=sqN+tN+i$ if $v_i=[a/c]_q$ is in position $i$ in $\M_3(q)$.  This ensures that, for any pair of integers $s,t=0,\ldots,p-1$, the set of vertices $w_l:l=sqN+tN+i,i=0,\ldots,N-1$ contains exactly one vertex in each of the fibres of the vertices $v_i$ of $\M_3(q)$.  Then, apart from the entries when $\Delta_{pq}\equiv0\pmod{pq}$, the matrix $A(pq)^2$ consists of $p^2\times p^2$ copies of $A(q)^2$. 
If $\Delta_{pq}\equiv0\pmod{pq}$, the $fg$ entry of $A(pq)^2-pqI$ is $0$.
The first row of $A(pq)^2-pqI$ consists of $p^2$ copies of  the first row of $A(q)^2$, apart from entries equal to $0$ corresponding to the vertices $[1+lq/0]_{pq}$ in positions $lpN$ for $l=0,\ldots,p-1$.  We can check that an automorphism takes $[1/0]_{pq}$ to $w_f$, and $[1+lq/0]_{pq}$ to a vertex in a position a multiple of $pN$ from $w_f$.
%%%%%%%%%%%%%%%%%%%%%%%%%%%%%%%%%%%%%%%%%%
We define the $N\times N$ matrices $B=A^2(q)$, $T=A^2(q)-qI,$ and the $pN\times pN$ matrix $D$, which consists of $p$ blocks of rows each comprising $p-1$ copies of $B$, with one copy of $T$ in the diagonal position.
\begin{eqnarray*}\label{pkmatrix}D=\begin{pmatrix}T&B&\cdots&B\\B&T&\cdots&B\\\vdots&\vdots&&\vdots\\B&B&\cdots&T\end{pmatrix}.\quad\text{ Then }
A^2(pq)-pqI=\begin{pmatrix}D&D&\cdots&D\\D&D&\cdots&D\\\vdots&\vdots&&\vdots\\D&D&\cdots&D\end{pmatrix}.
\end{eqnarray*} The rank of $D$ is less than or equal to $pN$. 
 So, from Lemma \ref{blocks}, $\sqrt{pq}$ or $-\sqrt{pq}$ or both  are eigenvalues of $\M(pq)$ with $\gamma(-\sqrt{pq})+\gamma(\sqrt{pq})\ge\frac{1}{2}(p^2N-pN)=\frac{1}{2}p(p-1)N$. 
\end{proof}
\noindent We can now prove Theorem \ref{B1}(ii).

\begin{proof}[Proof of Theorem \ref{B1}(ii)]

We put $p^{k-1}=q$.
The number of vertices of $\M(pq)$, and so the number of its eigenvalues,
 is $p^2N$.  Either $\sqrt{pq}$ or $-\sqrt{pq}$,
 or both, are eigenvalues of $\M_3(pq)$; $\gamma(-\sqrt{pq})+\gamma(\sqrt{pq})\ge p(p-1)N$,
 and  $\gamma(0)\ge(p-1)N$. 
 Also, from Lemma \ref{covereig}, $p\ \!\textnormal{sp}_3(q)\subset\textnormal{sp}_3(pq)$. The sum of the lower bound of the multiplicities of all the eigenvalues of $\M_3(pq)$ we have found is $p(p-1)N+(p-1)N+N=Np^2$, which is the number of eigenvalues of $\M_3(pq)$, so we take the lower bound for all multiplicities
 and there are no more eigenvalues. Since the entries on the main diagonal of $A(pq)$ are all zero, its trace is zero,
 so the sum of its eigenvalues  is 0.
Therefore, as the eigenvalues of ${\mathcal{M}}_3(q)$ also sum to 0, $\sqrt{pq}$ and 
$-\sqrt{pq}$ have the same multiplicity as eigenvalues of $\M(pq)$.  
\end{proof}

%%%%%%%%%%%%%%%%%%%%%%%%%%%%%%%%%%%%%%%%%%%%%%%
\section{The spectrum of $\M_3(n)$ for a composite $n$:  proof of Theorem 1.4}
%%%%%%%%%%%%%%%%%%%%%%%%%%%%%%%%%%%%%%%%%
To find the spectrum of $\M_3(n)$ for a composite $n$ we need the parallel product of maps introduced in Section \ref{3}.   From \cite[Theorem 2]{SW2002}, if the vertex valencies of two maps $\M_1$ and $\M_2$ are coprime, the underlying graph of their parallel product $\M$ is the tensor product $\G_1\times\G_2$ of their underlying graphs.  Then, from \cite[Theorem 4.2.12]{HJ}, $\textnormal{sp}(\M)=\textnormal{sp}(\M_1)\textnormal{sp}(\M_2)$.
\begin{lemma}\label{sqlm}If $l$ and $m$ are coprime integers, neither of which is twice an odd integer, then $\sqrt{lm}$ or $-\sqrt{lm}$ or both are eigenvalues of $\M_3(lm)$ with total multiplicity equal to half the number of its vertices.
\end{lemma}
\begin{proof}
From Theorem~\ref{A2}, ${\mathcal{M}}_3(lm)$ is a double covering of the parallel product of ${\mathcal{M}}_3(l)$ and $\M_3(m)$.  Let $u\in\Z, 1<u<lm$ be such that $u\equiv 1\!\!\pmod{l}$ and $u\equiv-1\!\!\pmod{m}$. The covering transformation takes both the vertices $w_f= [a/c]_{lm}$ and $w_{f'}=[ua/uc]_{lm}$of ${\mathcal{M}}_3(lm)$ to the vertex $v_i=([a/c]_l,[a/c]_m)$ of the parallel product of ${\mathcal{M}}_3(l)\text{ and }\M_3(m)$.  Let the number of vertices of ${\mathcal{M}}_3(lm)$ be $2V$.  Then the parallel product has $V$ vertices.
We order the vertices of $\M_3(lm)$ so that $f'=f+V$.  
 Let $w_g=[b/d]_{lm}$ also be a vertex of ${\mathcal{M}}_3(lm)$.  We compare the $fg$ and $f'g$ entries of the matrix $(A(lm))^2-lmI$.  We have $\Delta(fg)=ad-bc$ and $\Delta(f'g)=uad-buc$.  Then gcd$(\Delta(fg),lm)=\textnormal{gcd}(\Delta(f'g),lm)=r$. We note that $u\lambda(fg) u a+u\mu(fg) uc+\nu lm=1$, so we put $\lambda(f'g)=u\lambda(fg)$ and $\mu(f'g)=u\mu(fg)$, and obtain $\beta(f'g)=\beta(fg)$.  So the $A(lm)^2-lmI$ entries $fg$ and $f+V,g$ are equal. We can write $A^2(lm)-lmI$ as $2\times2$ blocks of an $V\times V$ matrix and apply Lemma \ref{blocks}.  So $\pm\sqrt{lm}$ are eigenvalues of $(A(lm))^2$ with $\gamma(-\sqrt{lm})+\gamma(\sqrt{lm})\ge V$.\end{proof}

It is now straightforward to prove Theorem~\ref{B2} 

\begin{proof}[Proof of Theorem \ref{B2}]

Both parts of this theorem now follow from Theorem \ref{A2}.  For part (ii), we also need Lemma~\ref{sqlm}, and we note that, as the trace of the adjacency matrix of the parallel product is zero, its eigenvalues sum to zero; therefore the eigenvalues $\sqrt{lm}$ and $-\sqrt{lm}$ have the same multiplicity since the eigenvalues of $\M_3(lm)$ also sum to zero.
\end{proof}
%%%%%%%%%%%%%%%%%%%%%%%%%%%%%%%%%%%%%%%%
\section{Ramanujan graphs}
A Ramanujan graph of degree $n$ is a regular graph for which 
 ${\lambda}_1<2\sqrt{n-1}$, where ${\lambda}_1$ is the graph eigenvalue with the second largest modulus.  They are the subject of much research as they make particularly good communication networks.
\begin{lemma}\label{t}
  If $p_1$ is the smallest prime dividing $n$, the eigenvalue ${\lambda}_1(n)$ of $\M_3(n)$ is, for  $n>4$, \[{\lambda}_1(n)=\begin{cases}\frac{1}{2}n&\text{if}\quad p_1=2,\\\frac{1}{3}n&\text{if}\quad p_1=3\\\frac{1}{\sqrt{p_1}}n&\text {otherwise}.\end{cases}\]
\end{lemma}
\begin{proof}In the spectrum of $p^k$, the largest eigenvalue is $p^k$, and, from Theorem \ref{B1}, the eigenvalue with the second largest modulus is $p^k/p^\frac{1}{2}$ for $p>3$, $p^k/2$ for $p=2$, and $p^k/3$ for $p=3$.
We proceed by induction on the number $r$ of prime divisors of $n$.  If $r=1, n=p_1^k$ and ${\lambda(n)}_1=n/\sqrt{p_1}$ or $n/2$ or $n/3$.  Now assume that the lemma is true for any integer with $r-1$ prime divisors.  Let $n$ be any integer with $r$ divisors the smallest of which is $p_1$.  Write $n=p_1^km$, so that $m$ is an integer with $r-1$ prime divisors the smallest of which is $p_2>p_1$.  Then, by Theorem~\ref{B2} and the induction assumption, the eigenvalue with the second largest modulus of  $\M_3(n)$ is $\max\{p^k_1(m/\sqrt{p_2}),({p_1^k}/{\sqrt{p_1}})m\}=n/\sqrt{p_1}$, which proves the result for $p_1>3$.  It is straightforward to prove in the same way the corresponding results for $p_1=2$ and $p_1=3$.
\end{proof}
From this lemma we recover a result proved in a very different way in \cite[Theorem 1]{DLM}, \cite[Theorem 4.2]{Gun}, and \cite[Theorem 1.12,(i)]{N}.
\begin{corollary}\label{R}The underlying graph of the Farey map $\M_3(n)$ is Ramanujan if and only if $n$ is prime, or equal to one of $4,6,8,9,10,12,14,15,21,27$ or $33$.
\end{corollary}
\begin{proof}
The map is Ramanujan if and only if ${\lambda}_1<2\sqrt{n-1}$.  Let $p_1$ be the smallest prime factor of $n$.  Then if $p_1>3$, that is equivalent to $n/\sqrt{p_1}<2\sqrt(n-1)$, which is true if and only if $n<4p_1$, that is if and only if $n=p_1$, so that $n$ is prime.
If $p_1=2,$ the map is Ramanujan if and only if $n/2<2\sqrt{n-1}$, which is equivalent to $n^2-16n+16<0$, which implies $n\le14$.  If $p_1=3$,   the map is Ramanujan if and only if $n/3<2\sqrt{n-1}$, which is equivalent to $n^2-36n+36<0$, which implies $n\le34$.
\end{proof}

%%%%%%%%%%%%%%%%%%%%%%%%%%%%%%%%%%%%%%%
\section{Maps defined by Hecke groups}
%%%%%%%%%%%%%%%%%%%%%%%%%%%%%%%%%%%%%
The Hecke group  $H^q$ is a discrete subgroup of infinite index in PSL$_2(\mathbb{Z}[{\lambda}_q])$ generated by the matrices $R=\begin{bmatrix} 0&-1\\1&0\end{bmatrix}$ and $S=\begin{bmatrix}1&\lambda_q\\0&1\end{bmatrix}$, where  ${\lambda}_q=2\cos\pi/q$.  
 The universal Hecke map $\widehat{\mathcal{M}}_q$ is the tessellation of the upper hyperbolic plane whose darts are the dart from infinity to zero and its images  under $H^q$. 
 If $q=3, \lambda_q=1, H^3$ is the modular group and $\widehat{\mathcal{M}}_3$ is the Farey tessellation.
 Hecke maps are quotients of $\widehat{\mathcal{M}}_q$ by the congruence subgroups of $H^q$.  We are particularly interested in the map $\mathcal{M}_q(n)$ of type $(q,n)$, which is the quotient of $\widehat{\mathcal{M}_q}$ by the subgroup $H^q(n)=H^q/(n)$ defined by the ideal $(n)=\{n(a+b{\lambda}_q):a,b\in\mathbb{Z}\}$ of $\mathbb{Z}[{\lambda}_q] $. 

As in \cite{IS,KaSi}, we will consider the maps ${\mathcal{M}}_4(n)$ and ${\mathcal{M}}_6(n)$ corresponding to the Hecke groups $H^4$ and $H^6$.  These are relatively straightforward to consider as ${\lambda}_4=\sqrt{2}$ and  ${\lambda}_6=\sqrt{3}$.  The faces of ${\mathcal{M}}_4(n)$ and ${\mathcal{M}}_6(n)$ are, respectively, quadrilaterals and hexagons.    As is shown in \cite{IS}, ${\mathcal{M}}_4(n)$ has two types of vertices:  even vertices, which are the equivalence classes of ordered pairs $\{(a',c'\sqrt{2}):(a',c')\in\Z\times\Z: \textnormal{gcd}(a',c',n)=1,\textnormal{gcd}(a,2,n)=1\;(a',c')\equiv \pm(a,c)\!\!\pmod{n}\}$, which we write $[a/c\sqrt{2}]_n$, and odd vertices, which are the equivalence classes of ordered pairs $\{(a'\sqrt{2},c'):(a',c')\in\Z\times\Z: \textnormal{gcd}(a',c',n)=1,\textnormal{gcd}(c',2,n)=1\;(a',c')\equiv \pm(a,c)\!\!\pmod{n}\}$, which we write $[a\sqrt{2}/c]_n$.
 If  $[a,c\sqrt{2}]_n$ and $[b\sqrt{2},d]_n$ with $a,b,c,d\in \mathbb{Z}/n\mathbb{Z}$ are two vertices of ${\mathcal{M}}_{4}(n)$, then those vertices are adjacent if and only if
 $ad-2bc\equiv \pm1\!\! \pmod n$. Odd vertices are adjacent to even vertices, and vice-versa.
 The vertex valency of  ${\mathcal{M}}_{4}(n)$ is $n$.
Replacing 2 by 3 and 4 by 6 gives analogous results for ${\mathcal{M}}_{6}(n)$.

  ${\mathcal{G}}_4(n)$ is the underlying graph of  ${\mathcal{M}}_4(n)$, and ${\mathcal{G}}_6(n)$ that of  ${\mathcal{M}}_6(n)$.
\begin{theorem}\label{Hecke} For odd $n$, ${\mathcal{G}}_4(n)$, is a double graph covering of ${\mathcal{G}}_3(n)$.  If $n$ is not a multiple of $3$, ${\mathcal{G}}_6(n)$ is a double graph covering of ${\mathcal{G}}_3(n)$. 
\end{theorem}
\begin{proof}
 Consider the mapping
\vspace{-2ex}\begin{eqnarray*}
\sigma:\text{ vertices of } {\mathcal{M}}_4(n)&\longrightarrow&\text{ vertices of } {\mathcal{M}}_3(n)\\
\text{  } [a/c\sqrt{2}]_n&\longmapsto&[a/c]_n\\
\text{  } [a\sqrt{2}/c]_n&\longmapsto&[a/c]_n
\end{eqnarray*}
Then $\sigma$ is a bijection between the odd vertices of ${\mathcal{M}}_4(n)$ and the vertices of ${\mathcal{M}}_3(n)$, and also between the even vertices of ${\mathcal{M}}_4(n)$ and the vertices of ${\mathcal{M}}_3(n)$. So it is a graph covering transformation between the underlying graphs, of order 2 as each vertex of ${\mathcal{M}}_3(n)$ has two pre-images.  
Replacing 2 by 3 gives the corresponding result for ${\mathcal{M}}_{6}(n)$ if 3 is not a factor of $n$.
\end{proof}
\noindent Figure~\ref{Heckemaps} shows the cube with skeleton ${\mathcal{G}}_4(3)$ as a double graph covering of the tetrahedron with skeleton ${\mathcal{G}}_3(3)$.  Both these graphs can be embedded as maps in the sphere.  This is not a map covering as there is no mapping from the $6$ faces of ${\mathcal{M}}_4(3)$ to the $4$ faces of ${\mathcal{M}}_3(3)$. 

\begin{figure}[h]
\begin{center}
\begin{tikzpicture}[scale=0.9]
%M_4(3)%
\draw (0,0) rectangle (2.5,2.5);
\draw(0,2.5)--(1.5,3)--(4,3)--(2.5,2.5);
\draw((4,3)--(4,0.5)--(2.5,0);
\draw[dashed](0,0)--(1.5,0.5)--(4,0.5);
\draw[dashed](1.5,0.5)--(1.5,3);
\draw [fill=white] (0,0) circle [radius=0.12];
\node[ below left] at  (0,0){$[1/0\sqrt{2}]_3$};
\draw [fill] (2.5,0) circle [radius=0.1];
\node[below right] at  (2.5,0){$[0\sqrt{2}/1]_3$};
\draw [fill=white](2.5,2.5) circle [radius=0.12];
\node[ below right] at  (2.5,2.5){$[2/\sqrt{2}]_3$};
\draw [fill](0,2.5) circle [radius=0.1];
\node[below left] at  (0,2.5){$[\sqrt{2}/1]_3$};
\draw [fill=white] (4,0.5) circle [radius=0.12];
\node[ above right] at  (4,0.5){$[1/\sqrt{2}]_3$};
\draw [fill] (1.5,0.5) circle [radius=0.1];
\node[above left] at  (1.65,0.5){$[2\sqrt{2}/1]_3$};
\draw [fill=white](1.5,3) circle [radius=0.12];
\node[ above left] at  (1.5,3){$[0/\sqrt{2}]_3$};
\draw [fill] (4,3) circle [radius=0.1];
\node[above right] at  (4,3){$[\sqrt{2}/0]_3$};

\draw [fill=white](-4.75,1.7) circle [radius=0.12];
\node[ right] at  (-4.65,1.7){even vertex};

\draw [fill](-4.75,1) circle [radius=0.12];
\node[ right] at  (-4.65,1){odd vertex};

%$\M_3(3)$

\draw(7,0)--(8.25,3)--(10,0.5)--(9.5,0)--(7,0);
\draw(8.25,3)--(9.5,0);
\draw[dashed](7,0)--(10,0.5);
\draw [fill=lightgray](7,0) circle [radius=0.12];
\node[below left] at  (7,0){$[1/0]_3$};
\draw [fill=lightgray](8.25,3) circle [radius=0.12];
\node[above right] at  (8.25,3){$[2/1]_3$};
\draw [fill=lightgray](10,0.5) circle [radius=0.12];
\node[right] at  (10,0.5){$[1/1]_3$};
\draw [fill=lightgray](9.5,0) circle [radius=0.12];
\node[right] at  (9.5,0){$[0/1]_3$};

\end{tikzpicture}
\caption{The Hecke map $\M_4(3)$ and the Farey map $\M_3(3)$.}
\label{Heckemaps}
\end{center}
\end{figure}
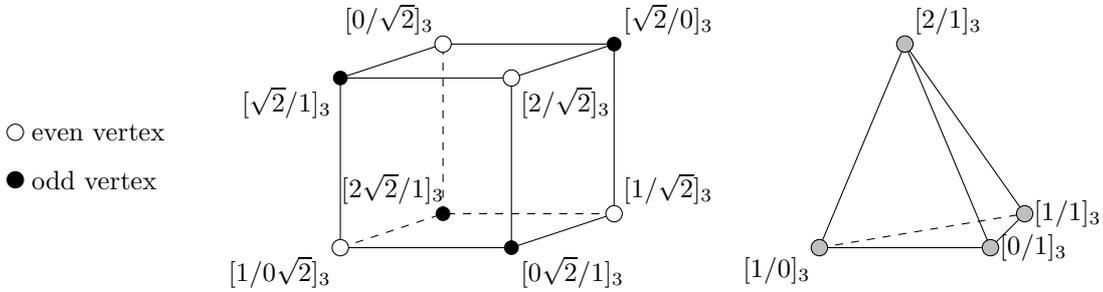
\begin{corollary} The spectrum of ${\mathcal{M}}_4(n)$ for odd $n$,
 and that of ${\mathcal{M}}_6(n)$ for $3\nmid n$, is the multiset $-\textnormal{sp}_3(n)\cup\textnormal{sp}_3(n)$, 
 where $\textnormal{sp}_3(n)$ is the spectrum of ${\mathcal{M}}_3(n)$.  
\end{corollary}
\begin{proof}
From Theorem~\ref{Hecke} and from Lemma \ref{covereig}, if ${\lambda}_i$ is an eigenvalue of ${\mathcal{M}}_3(n)$, then it is also an eigenvalue of ${\mathcal{M}}_4(n)$.  Moreover the vertices adjacent to an even vertex are all odd, and vice versa, so that ${\mathcal{M}}_4(n)$ is bijective.  Then from Theorem 8.8.2 of \cite{GR}, $-{\lambda}_i$ also an eigenvalue of ${\mathcal{M}}_4(n)$.  Since ${\mathcal{M}}_4(n)$ has twice as many vertices as ${\mathcal{M}}_3(n)$ it has twice as many eigenvalues, so there are no more eigenvalues and the result follows.  
Replacing 4 by 6, we have the corresponding result for ${\mathcal{M}}_{6}(n)$ if 3 is not a factor of $n$.
\end{proof}
The following corollary now follows immediately from the result for ${\mathcal{M}}_3(p)$.
\begin{corollary} The underlying graphs of the Hecke maps ${\mathcal{M}}_{4}(p)$ and ${\mathcal{M}}_{6}(p)$ for a prime $p$ are Ramanujan.
\end{corollary}
We also obtain the diameter of certain Hecke maps, partially recovering a result in \cite{KaSi}.  We denote the shortest distance between 2 vertices $u$ and $v$ of a graph or map as $d(u,v)$.  The \emph{diameter} is the largest value of $d(u,v)$.
\begin{theorem} Both ${\mathcal{M}}_4(n)$ (for odd $n$) and ${\mathcal{M}}_6(n)$  (for $3\nmid n$) have diameter $4$.\end{theorem}
\begin{proof} Define $\sigma$ as in Theorem~\ref{Hecke}, and let $P=\langle w_1,w_2,\ldots,w_k\rangle$ be a path in $\M_3(n)$.  
Then $\sigma^{-1}(w_i)$  consists of an even vertex of $\M_4(n)$, $v_i$, and an odd vertex $u_i$.  As even vertices are adjacent to odd vertices and vice-versa, $P$ is lifted by ${\sigma}^{-1}$ to two paths in $\M_4(n)$:  $P_1=\langle v_1,u_2,v_3,\ldots\rangle$, and $P_2=\langle u_1,v_2,u_3\ldots\rangle.$
 A lift of a path of even length in $\M_3(n)$ will join two vertices of the same  parity in $\M_4(n)$, and a lift of a path of odd length two vertices of opposite parity. 
 Now let $v_a$ and $v_b$ be any two distinct vertices of $\M_4(n)$, and let  $w_a=\sigma(v_a)$ and $w_b=\sigma(v_b)$. 
If $w_a=w_b$, $v_a$ and $v_b$ are of opposite parities.
  There is a path of length 3 round the edges of a triangle incident to $w_a$; this lifts to a path of length 3 between $v_a$ and $v_b$, so $d(v_a,v_b)\le3$. 
   From \cite{SS}, the diameter of $\M_3(n)$ for any $n$ is 3, so if $w_a\neq w_b$, $d(w_a,w_b)\le3$.
  If $P$ is a path of shortest length between $w_1$ and $w_i$, the vertices preceding $w_{i-1}$ cannot be adjacent to $w_i$. 
 The edge $w_{i-1}w_i$ is incident to two triangles.  Let $w$ be the third vertex of one of these triangles. 
 Then there is a path $\langle w_1,\ldots,w_{i-1},w,w_i\rangle$ from $w_1$ to $w_i$ of length $d(w_1,w_i)+1$. 
 So since $d(w_a,w_b)\le3$, there is always both a path of odd length less than $4$ and a path of even length less than or equal to $4$ between  $w_a$ and $w_b$, which lifts to a path of length less than or equal to $4$ between $v_a$ and $v_b$  whether or not they are of the same parity.
The result for $\M_6(n)$ if $3\nmid n$ follows similarly. 
\end{proof}
The further result that this is also true for all $n>6$ is given by \cite[Theorem 12]{KaSi}. 
\section{acknowledgments}
The author would like to thank Ian Short for his invaluable support and guidance, Nick Gill and Jozef \v{S}ir\'{a}\v{n} for helpful discussions, and the anonymous referee for their comments, which have greatly improved this paper.

%%%%%%%%%%%%%%%%%%%%%%%%%%%%%%%%%%%%%%%%%%%%%%%%%%%%%%%%%%%%%%%

%%%%%%%%%%%%%%

\begin{thebibliography}{99}
\bibitem{A} T.M. Apostol, {\it Introduction to Analytic Number Theory}, Undergraduate Texts in Mathematics, Springer, New York, 1976.
\bibitem{BrNe}A. Breda D'Azevedo\ and\ R. Nedela, Join and intersection of hypermaps, Acta Univ. M. Belii Ser. Math. No. 9 (2001), 13--28.
\bibitem{DFMRS} C. Dalf\'{o}, M. A. Fiol,\ M. Miller, J. Ryan and\ J. \v{S}ir\'{a}\v{n}, An algebraic approach to lifts of digraphs, Discrete Appl. Math. {\bf 269} (2019), 68--76.
\bibitem{DFPS}C. Dalf\'{o}, M. A. Fiol,\ S. Pavl\'{i}kov\'{a} and\ J. \v{S}ir\'{a}\v{n}, Spectra and eigenspaces of arbitrary lifts of graphs, Acta Math. Univ. Comenian. (N.S.) {\bf 88} (2019), no.~3, 593--600. 
\bibitem{DFS} C. Dalf\'{o}, M.A. Fiol\ and\ J. \v{S}ir\'{a}\v{n}, The spectra of lifted digraphs, J. Algebraic Combin. {\bf 50} (2019), no.~4, 419--426. 
\bibitem{DLM} M. DeDeo, D. Lanphier\ and\ M. Minei, The spectrum of Platonic graphs over finite fields, Discrete Math. {\bf 307} (2007), no.~9-10, 1074--1081. 
\bibitem{GR}C. Godsil\ and\ G. Royle, {\it Algebraic graph theory}, Graduate Texts in Mathematics, 207, Springer-Verlag, New York, 2001.
\bibitem{Gun} P.E. Gunnells, Some elementary Ramanujan graphs, Geom. Dedicata {\bf 112} (2005), 51--63.
\bibitem{HJ} R.A. Horn\ and\ C.R. Johnson, {\it Topics in matrix analysis}, Cambridge University Press, Cambridge, 1991.
\bibitem{I19}I. Ivrissimtzis, N. Peyerimhoff\ and\ A. Vdovina, Trivalent expanders, $(\Delta-Y)$-transformation, and hyperbolic surfaces, Groups Geom. Dyn. {\bf 13} (2019), no.~3, 1103--1131.
\bibitem{IS}	I. Ivrissimtzis\ and\ D. Singerman, Regular maps and principal congruence subgroups of Hecke groups, European J. Combin. {\bf 26} (2005), no.~3-4, 437--456.
\bibitem{Hu} K. Hu, R. Nedela\ and\ N-E. Wang, Branched cyclic regular coverings over platonic maps, European J. Combin. {\bf 36} (2014), 531--549.
\bibitem{Hu16} K. Hu\ et al., Non-abelian almost totally branched coverings over the platonic maps, European J. Combin. {\bf 51} (2016), 1--11

\bibitem{TMOS} G.A. Jones\ and\ D. Singerman, Theory of maps on orientable surfaces, Proc. London Math. Soc. (3) {\bf 37} (1978), no.~2, 273--307.
\bibitem{KaSi} D. Kattan\ and D. Singerman, The diameter of some Hecke-Farey maps, Albanian Journal of Mathematics {\bf 15}(2021), no.~1, 39-60.
\bibitem{MMW}P. McMullen, B. Monson\ and\ A.I. Weiss, Regular maps constructed from linear groups, European J. Combin. {\bf 14} (1993), no.~6, 541--552
\bibitem{N}B. Nica, Unimodular graphs and Eisenstein sums. J. Algebraic Combin. 45 (2017), no. 2, 423--454.
\bibitem{UT} D. Singerman, Universal tessellations, Rev. Mat. Univ. Complut. Madrid {\bf 1} (1988), no.~1-3, 111--123.
\bibitem{SS} D. Singerman\ and\ J. Strudwick, The Farey maps modulo $n$, Acta Math. Univ. Comenian. (N.S.) {\bf 89} (2020), no.~1, 39--52.
\bibitem{SS16} D. Singerman\ and\ J. Strudwick, Petrie polygons, Fibonacci sequences and Farey maps, Ars Math. Contemp. {\bf 10} (2016), no.~2, 349--357. 
\bibitem{SW2002} S.E. Wilson, Families of regular graphs in regular maps, J. Combin. Theory Ser. B {\bf 85} (2002), no.~2, 269--289.
\bibitem{SW1994} S.E. Wilson, Parallel products in groups and maps, J. Algebra {\bf 167} (1994), no.~3, 539--546.
\end{thebibliography}
\end{document}